\documentclass[a4paper,twoside,12pt]{article}
\usepackage[english]{babel}
\usepackage{graphicx}
\usepackage{bbding}
\usepackage[T1]{fontenc}
\usepackage{kpfonts}
\usepackage{fancyhdr}
\usepackage{amsmath,amssymb,color,epsfig}
\usepackage{mathrsfs}
\usepackage{dsfont}
\usepackage{upgreek}
\usepackage{fancyhdr}

\usepackage{enumitem}
\newtheorem{theorem}{Theorem}[section]
\newtheorem{proposition}[theorem]{Proposition}
\newtheorem{corollary}[theorem]{Corollary}

\newtheorem{lemma}[theorem]{Lemma}

\def\R{{\Bbb R}}

\def\1{\mathds{1}}

\title{A harmonic mean inequality for the $q-$gamma and $q-$digamma functions}
\author{Mohamed BOUALI}
\date{}
\begin{document}
\maketitle
\begin{abstract} We prove amongs others results that the harmonic mean of $\Gamma_q(x)$ and $\Gamma_q(1/x)$ is greater than or equal to $1$
for arbitrary $x > 0$ and  $q\in J$ where $J$ is a subset of $[0,+\infty)$. Also, we prove that for there is $p_0\in(1,9/2)$, such that for $q\in(0,p_0)$, $\psi_q(1)$ is the minimum of the harmonic mean of $\psi_q(x)$ and  $\psi_q(1/x)$ for $x > 0$ and for $q\in(p_0,+\infty)$, $\psi_q(1)$ is the maximum. Our results generalize some known inequalities due to Alzer and Gautschi.
\end{abstract}
\section{Preliminaries and notations}
The $q$-analogue of the $\Gamma$ function, denoted by $\Gamma_q(x)$, was introduced by \cite{jak} as
\begin{equation}\label{e1}\Gamma_q(x)=(1-q)^{1-x}\prod_{n=0}^\infty\frac{1-q^{n+1}}{1-q^{n+x}},\qquad 0<q<1,\end{equation}
and
\begin{equation}\label{e2}\Gamma_q(x)=(q-1)^{1-x}q^{x(x-1)/2}\prod_{n=0}^\infty\frac{1-q^{-(n+1)}}{1-q^{-(n+x)}},\qquad q>1,\end{equation}
for $x>0$. It was proved in \cite{mok} that $\Gamma(x)=\lim_{q\to 1^-}\Gamma_q(x)=\lim_{q\to 1^+}\Gamma_q(x).$

Similarly the $q$-digamma or $q$-psi function $\psi_q(x)$ is defined by $\psi_q(x)=\frac{\Gamma_q'(x)}{\Gamma_q(x)}.$ The derivatives $\psi'_q, \psi_q''$,...are  called  the $q$-polygamma functions. In \cite{krat}, it was shown that $\lim_{q\to 1^-}\psi_q(x)=\lim_{q\to1+}\psi_q(x)=\psi(x)$. From the definitions \eqref{e1} and \eqref{e2} one can easily deduce that
 for $0<q<1$, \begin{equation}\label{e3}\psi_q(x)=-\log(1-q)+(\log q)\sum_{n=1}^\infty\frac{q^{nx}}{1-q^n}=-\log(1-q)+(\log q)\sum_{n=1}^\infty\frac{q^{n+x}}{1-q^{n+x}},\end{equation}
and for
$q>1$, \begin{equation}\label{e4}\psi_q(x)=-\log(q-1)+(\log q)\Big[x-\frac12-\sum_{n=1}^\infty\frac{q^{-nx}}{1-q^{-n}}\Big]=-\log(q-1)+(\log q)\Big[x-\frac12-\sum_{n=1}^\infty\frac{q^{-(n+x)}}{1-q^{-(n+x)}}\Big],\end{equation}
Differentiation of \eqref{e3} and \eqref{e4} gives
\begin{equation}\label{e5}\psi_q'(x)=\left\{\begin{aligned}&(\log q)^2\sum_{n=1}^\infty\frac{n q^{nx}}{1-q^n},\qquad 0<q<1\\
&\log q+(\log q)^2\sum_{n=1}^\infty\frac{n q^{-nx}}{1-q^{-n}},\qquad q>1,\end{aligned}\right.\end{equation}
or equivalently
\begin{equation}\label{e5}\psi_q'(x)=\left\{\begin{aligned}&(\log q)^2\sum_{n=1}^\infty\frac{ q^{n+x}}{(1-q^{n+x})^2},\qquad 0<q<1\\
&\log q+(\log q)^2\sum_{n=1}^\infty\frac{q^{-(n+x)}}{(1-q^{-(n+x)})^2},\qquad q>1,\end{aligned}\right.\end{equation}
One deduces that $\psi_q$ is strictly increasing on $(0,+\infty)$.

A short calculation gives for $x>0$ and $q>0$
 \begin{equation}\label{e6}\Gamma_q(x)=q^{(x-1)(x-2)/2}\Gamma_{\frac1q}(x).\end{equation}
If we take logarithm of both sides of \eqref{e6} and then differentiate, we find
$$\psi_q(x)=(x-\frac32)\log q+\psi_{\frac1q}(x).$$

In \cite{koc}, it was proved that the function $q\mapsto \psi_q(1)$ decreases on $(0,1)$, and $q\mapsto \psi_q(2)$ increases on $(0,1)$. Then $\psi_q(1)\leq\psi_0(1)=0$ and $0=\psi_0(2)\leq\psi_q(2)$.
\begin{proposition} The function $\psi_q$ $(0<q)$ has a uniquely determined positive zero on $(1,\frac32)$, which we denote by $x_q$.
\end{proposition}
{\bf Proof.} Let $q\in(0,1)$, the function $\psi_q$ is strictly increasing on $(0,+\infty)$ and $\psi_q(1)\leq 0$ and $\psi_q(2)\geq 0$. Then $\psi_q$ has a unique zero in $(1,2)$. Furthermore,
 in \cite{bat5} Corollary 2.3, it was shown that $\log(\frac{q^{x+\frac12}-1}{q-1})<\psi_q(x+1)$ for all $x>0$ and $q>0$. Then for $x=\frac12$ and $q\in(0,1]$ we get
$\psi_q(\frac32)>0$. Moreover, for $q> 1$, $\psi_q(\frac32)=\psi_{\frac1q}(\frac32)$. Since, the function $\psi_q$ is strictly monotone on $(0,+\infty)$ for all $q>0$, hence $x_q\in(1,\frac32)$.

In \cite{alz2} H. Alzer proved for $x>0$, \begin{equation}\label{ik}\psi''(x)+(\psi'(x))^2>0.\end{equation} The author rediscovered it in \cite{bat1} and used it to prove interesting inequalities for the digamma function, see \cite{alz1, alz2, alz3}. Alzer and Grinshpan in \cite{alz3} obtained a $q$-analogue of \eqref{i} and proved that, for $q>1$ and all $x>0$,
\begin{equation}\label{i}\psi_q''(x)+(\psi'_q(x))^2>0.\end{equation}
The  author  in \cite{bat4} provided  another $q$-extension  of \eqref{i} and proved that
\begin{equation}\label{i}\psi_q''(x)+(\psi'_q(x))^2-\log q (\psi'_q(x))>0,\end{equation}
for all $q>0$ and all $x>0$. The following lemma is due to Alzer \cite{alz3}.
\begin{lemma}\label{l1}
For every $q>0$ and $x\geq 1$,
$$x\psi'_q(x)+2\psi_q(x)\geq 0.$$
\end{lemma}

We start with the following lemma which proves the convergence as $q\to 1$ of the $mth$ derivatives $\psi_q^{(m)}(x)$ to $\psi^{(m)}(x)$ for all $x>0$.
\begin{lemma}\label{der} The function $\psi_q^{(m)}(x)$ converges uniformly to $\psi^{(m)}(x)$ as $q\to 1$ on every compact of $(0,+\infty)$  for all $m\geq 0$, where $\psi_q^{(m)}$ respectively $\psi^{(m)}$ is the $mth$ derivatives of the psifunction respectively the $q$-psifunction.
\end{lemma}
{\bf Proof.} Let $h(u)=u/(1-u)$, then for every $m\geq 1$, and $u\neq 1$, $$h^{(m)}(u)=\frac{(-1)^{m+1}m!}{(1-u)^{m+1}}.$$
Firstly, we prove the lemma for $q\in(0,1)$. We have
$$\psi_q(x)=-\log(1-q)+(\log q)\sum_{n=0}^\infty h(q^{n+x}).$$
Recall the Fa\`a di Bruno formula for the $nth$ derivative of $f\circ g$, see for instance \cite{faa}. For $m\geq 0$,
$$(f\circ g)^{(m)}(t)=\sum_{k=0}^m\frac{m!}{k_1!k_2! ...k_m!}f^{(k)}(g(t))\prod_{i=1}^m(\frac{g^{(i)}(t)}{i!})^{k_i},$$
where $k=k_1+...+k_m$ and summation is over all naturel integers $k_1,...,k_m$\newline such that $k_1+2k_2+...+mk_m=m$. In our case,
$g^{(i)}(x)=(\log q)^iq^{x+n}$ and $f^{(k)}(x)=h^{(k)}(x)$. Then,
$$\big(h(q^{n+x})\big)^{(m)}=(\log q)^m\sum_k\frac{m!(-1)^{k+1}k!}{\prod_{i=1}^mk_i!\prod_{i=1}^m(i!)^{k_i}}\frac{q^{k(n+x)}}{(1-q^{n+x})^{k+1}}.$$
Remark that  for $k\leq m-1$, $x>0$  and $n\geq 0$, $\lim_{q\to 1^-}(\log q)^{m+1}\frac{q^{k(n+x)}}{(1-q^{n+x})^{k+1}}=0$. That is to say, the only term which contributes in the limit of  $(\log q)^{m+1}\frac{q^{k(n+x)}}{(1-q^{n+x})^{k+1}}$ as $q\to 1^-$ is when $k=m$. In this case, we have $k_1=m$, $k_2=k_3=...=k_{m}=0$. Then for $x>0$, $n\geq 0$ and $m\geq 1$, we get $$\lim_{q\to 1^-}(\log q)\big(h(q^{n+x})\big)^{(m)}=\lim_{q\to 1^-}m!\frac{(\log q)^{m+1}q^{m(n+x)}}{(1-q^{n+x})^{m+1}}=\frac{(-1)^{m+1}m!}{(n+x)^{m+1}}.$$
One shows that the functions $q\mapsto-q^y\log(q^y)/1-q^y$ and $y\mapsto-\log(q^y)/1-q^y$ increases and the functions $y\mapsto -q^y\log(q^y)/1-q^y$ and $q\mapsto-\log(q^y)/1-q^y$ decreases for all $y>0$ and $q\in(0,1)$.

Let $a,b>0$, then for all $x\in [a,b]$, $n\geq 0$ and all $q\in[1/2,1]$, we have
$$\big|(\log q)^{m+1}\frac{q^{k(n+x)}}{(1-q^{n+x})^{k+1}}\big|\leq \frac{(\log 2)^{n+b}}{1-2^{-n-b}}\frac{(\log 2)^{m-k}}{(n+a)^{k+1}},$$
Since, $m\geq 1$ then $k\geq 1$ and
$$\big|(\log q)\big(h(q^{n+x})\big)^{(m)}\big|\leq C_m \frac{(\log 2)^{n+b}}{1-2^{-n-b}}\frac{1}{(n+a)^2},$$
where $C_m$ is some constant independent of $n$.
This implies $\displaystyle\sum_{n=0}^\infty(\log q)\big(h(q^{n+x})\big)^{(m)}$ converges uniformly for  $(x,q)\in[a,b]\times[1/2,1]$ and $m\geq1$,  moreover $$\lim_{q\to 1^-}\sum_{n=0}^\infty(\log q)\big(h(q^{n+x})\big)^{(m)}= \sum_{n=0}^\infty\frac{(-1)^{m+1}m!}{(n+x)^{m+1}}=\psi^{(m)}(x).$$
Assume $q>1$. We saw that $\psi_q(x)=(x-3/2)\log q+\psi_{\frac1q}(x)$, and
$\psi^{(m)}_q(x)=\log q+\psi^{(m)}_{\frac1q}(x)$ if $m=1$,  $\psi^{(m)}_q(x)=\psi^{(m)}_{\frac1q}(x)$ if $m\geq 2$. The result follows from the case $q\in(1/2,1)$.

The case $m=0$ is proved in \cite{krat}.

\begin{lemma}\label{l2} For every $q>0$, the functions $f(x)=x\psi_q(x)$ and $g(x)=x\Gamma'_q(x)$ increase on $[1,+\infty)$.
\end{lemma}
{\bf Proof.} By differentiation we find $f'(x)=\psi_q(x)+x\psi'_q(x)$. By Lemma \ref{l1}, we get, for every $x\geq 1$, $f'(x)\geq-\psi_q(x)$. If $x\in[1,x_q]$ then $f'(x)\geq 0$ and then $f$ increases. For $x\geq x_q$, the functions $x\mapsto \psi_q(x)$ and $x\mapsto x$ are positive and increase then $f$ increases too.

To study the monotony of $g$, firstly remark that, $g(x)=x\psi_q(x)\Gamma_q(x)$. From the previous item and for $1\leq x\leq x_q$ we have the function $-x\psi_q(x)$ is positive and decreases and the function $\Gamma_q(x)$ is positive and decreases. Then $g$ increases too.

If $x\geq x_q$, since $g(x)=x\Gamma_q'(x)$ which is a product of two positive increasing functions then $g$ is an increasing function.
\begin{lemma}\label{l3} The function $x\psi'_q(x)$ is strictly decreasing on $(0,+\infty)$ for all $q\in(0,1)$, and the function $x^2\psi'_q(x)$ increases on $(0,+\infty)$ for all $q>1$.
\end{lemma}
{\bf Proof.}
1) Let $u(x)=x\psi_q'(x)$. By the integral representation of the $q$-digamma function we have,
$$u'(x)=\psi_q'(x)+x\psi_q''(x)=\int_0^\infty(1-t x)e^{-xt}\frac{t}{1-e^{-t}}d\gamma_q(t).$$
The function $t\mapsto t/(1-e^{-t})$ increases on $(0,+\infty)$. Then, for $t\in(0,\frac1x)$,
$$(1-t x)e^{-xt}\frac{t}{1-e^{-t}}\leq (1-t x)e^{-xt}\frac{1}{x(1-e^{-\frac 1x})}.$$
For $t\in(\frac 1x,+\infty)$, $1-tx\leq 0$ and $1/x(1-e^{-\frac1x})\leq t/(1-e^{-t})$. Hence
$$\int_0^\infty(1-t x)e^{-xt}\frac{t}{1-e^{-t}}d\gamma_q(t)\leq \frac{1}{x(1-e^{-\frac 1x})}\int_0^\infty(1-t x)e^{-xt}d\gamma_q(t).$$
Since,
$$-\log q\frac{q^x}{1-q^x}=\int_0^\infty e^{-x t}d\gamma_q(t).$$
$$-(\log q)^2\frac{xq^x}{(1-q^x)^2}=\int_0^\infty -xte^{-x t}d\gamma_q(t).$$
Then,
$$\int_0^\infty(1-t x)e^{-xt}\frac{t}{1-e^{-t}}d\gamma_q(t)\leq-\log q\frac{\frac1x}{1-e^{-\frac1x}}\frac{q^x}{1-q^x}\frac{1+x\log q-q^x}{1-q^x}<0,$$
for all $x\in(0,+\infty)$ and $q\in(0,1)$. Then, $u$ is strictly decreasing on $(0,+\infty)$.

2) Assume $q>1$, and let $$\varphi(x)=x^2\psi'_q(x),$$
then, $$\varphi'(x)=x(2\psi'_q(x)+x\psi''_q(x))\geq x(2\psi'_q(x)-x(\psi'_q(x))^2+x(\log q)\psi'_q(x)).$$
Thus,
$$\varphi'(x)\geq x\psi'_q(x)(2-x\psi'_q(x)+(\log q)x)= x\psi'_q(x)(2-x\psi'_{\frac1q}(x)).$$
We saw that for $q>1$, the function $x\mapsto x\psi'_{\frac1q}(x)$ decreases on $(0,+\infty)$, then for every $x\geq 1$,
$$\varphi'(x)\geq   x\psi'_q(x)(2-\psi'_{\frac1q}(1)).$$ On the other hand the function $\alpha\mapsto\psi'_\alpha(1)$ increases on $(0,1)$ and by Lemma \ref{der}, $\lim_{\alpha\to 1}\psi'_\alpha(1)=\psi'_1(1)=\frac{\pi^2}6$. Hence,
$$\varphi'(x)\geq   x\psi'_q(x)(2-\frac{\pi^2}6)\geq 0.$$

\begin{proposition}\label{prou} The functions $G_q(x)=\frac{\psi'_q(x)}{\psi_q(x)}$ and $\varphi_q(x)=x\frac{\psi'_q(x)}{\psi_q(x)}$ decrease on $(1,+\infty)$ for all $q\in(0,1]$.
\end{proposition}
{\bf Proof.} By differentiation we get
 $$G'_q(x)=\frac{\psi''_q(x)\psi_q(x)-(\psi'_q(x))^2}{(\psi_q(x))^2}.$$
 If $x\in(x_q,+\infty)$ then $\psi_q(x)\geq 0$ and $\psi''_q(x)\leq 0$ and the result follows. That is to say, the function $G_q(x)=\frac{\psi'_q(x)}{\psi_q(x)}$ decreases on $(x_q,+\infty)$ for all $q\in(0,1)$.

In order to prove $G'_q(x)\leq 0$ for $x\in(1,x_q)$, it is enough to show that for all $q\in(0,1)$, $$\psi''_q(x)\psi_q(x)-(\psi'_q(x))^2\leq0.$$
 For $x>0$ and $q>0$, the inequality $\psi''_q(x)\geq -(\psi'_q(x))^2+(\log q)\psi'_q(x)$ is proved in \cite{bat3}. So,
$$\psi''_q(x)\psi_q(x)-(\psi'_q(x))^2\leq -(\psi'_q(x))^2\psi_q(x)+(\log q)\psi'_q(x)\psi_q(x)-(\psi'_q(x))^2.$$
Which is equivalent to
$$\psi''_q(x)\psi_q(x)-(\psi'_q(x))^2\leq -\psi'_q(x)\Big(\psi'_q(x)\psi_q(x)+\psi'_q(x)-(\log q)\psi_q(x)\Big).$$
For $x\in(1,x_q)$ and $q\in(0,1)$, let's define $$\theta_q(x)=\psi'_q(x)\psi_q(x)+\psi'_q(x)-(\log q)\psi_q(x).$$
Differentiation of $\theta_q(x)$ gives
$$\theta'_q(x)=\psi''_q(x)\psi_q(x)+\psi''_q(x)+(\psi_q'(x))^2-(\log q)\psi'_q(x).$$
Hence, $\theta'_q(x)\geq\psi''_q(x)\psi_q(x)\geq 0$ for all $x\in(1,x_q)$.
Thus, $\theta_q(x)$ increases on $(1,x_q)$, and
$$\theta_q(x)\geq \theta_q(1)=\psi'_q(1)\psi_q(1)+\psi'_q(1)-(\log q)\psi_q(1).$$
It remains to show that the right hand side is positive.
We have, $$\theta_q(1)=(\psi'_q(1)-\log q)\psi_q(1)+\psi'_q(1).$$
In \cite{bat5}, it is proved that the function $F_q(x)=\psi_q(x+1)-\log(\frac{1-q^{x+\frac12}}{1-q})$ is completely monotonic on $(0,+\infty)$ for all $q>0$. Then, $F'_q(x)\leq 0$, which implies that $$\psi'_q(x+1)\leq -(\log q)\frac{q^{x+\frac12}}{1-q^{x+\frac12}}.$$
Then, $$\theta_q(1)\geq (-(\log q)\frac{q^{\frac12}}{1-q^{\frac12}}-\log q)\psi_q(1)+\psi'_q(1).$$
On other words,
 $$\theta_q(1)\geq \frac{\log q}{\sqrt q-1}\Big(\psi_q(1)+\frac{\sqrt q-1}{\log q}\psi'_q(1)\Big).$$
 Since, $ \frac{\log q}{\sqrt q-1}\geq 0$ for all $q\in(0,1)$. It is enough to prove that for all $q\in(0,1)$, $$u(q):=\psi_q(1)+\frac{\sqrt q-1}{\log q}\psi'_q(1)\geq 0.$$
 Using the series expansion of the functions $\psi_q(1)$ and $\psi'_q(1)$, we get
 $$u(q)=\sum_{n=1}^\infty\frac{q^n}{n(1-q^n)}\Big(1-q^n+n\log q+n^2(\sqrt q-1)\log q\Big).$$
 For $x\geq 1$, and $q\in(0,1)$, defines the function
 $$g(x)=1-q^x+x\log q+x^2(\sqrt q-1)\log q.$$
 By differentiation we find
 $$g'(x)=\left(1-q^x-2 x+2 x\sqrt{q} \right) \log q,$$
 $$g''(x)=\left(-2+2 \sqrt{q}-q^x \log q\right)\log q,$$
 $$g'''(x)=-q^x(\log q)^3.$$
 For $q\in(0,1)$, we have $g'''(x)>0$ for all $x\geq 1$, then
  $$g''(x)\geq g''(1)=-\left(2-2 \sqrt{q}+q \log q\right)\log q.$$
  Let $v(q)=2-2 \sqrt{q}+q \log q$, then $v'(q)=-\frac1{\sqrt q}+1+\log q$ and $v''(q)=\frac1q+\frac12q^{-\frac32}\geq 0$ and $v'(q)\leq v'(1)=0$. Thus, $v(q)$ decreases and $v(q)\geq v(1)=0$. Which implies that, $g''(x)\geq 0$, and for $q\in (0,1)$, $$g'(x)\geq g'(1)=-(\sqrt q-1)^2 \log q>0,$$
  hence, $$g(x)\geq g(1)=1-q+\sqrt q\log q.$$
  Setting $w(q)=1-q+\sqrt q\log q.$ Then,
  $$w'(q)=-1+\frac1{\sqrt q}+\frac1{2\sqrt q}\log q,$$ and
  $$w''(q)=-\frac14 q^{-\frac32}\log q>0.$$
  Then, $w'(q)\leq w'(1)=0$, and $w(q)\geq w(1)=0$.
  This implies $g(n)\geq 0$ for all $n\in\Bbb N$, and $u(q)\geq 0$ for all $q\in(0,1)$. On other words, $\theta_q(1)\geq 0$ for all $q\in(0,1)$.
Which gives the desired result.

  2) To prove that $\varphi_q(x)=x\frac{\psi'_q(x)}{\psi_q(x)}$ decreases. Remark that for $x\in(1,x_q)$ the function $-G_q(x)$ increases and is positive. Then $-\varphi_q(x)$ increases on $(1,x_q)$.

  For $x\in (x_q,+\infty)$ and $q\in(0,1)$ and by Lemma \ref{l3} we have $x\psi'_q(x)$ decreases and is positive. Moreover, $\frac1{\psi_q(x)}$ decreases and is positive for all $q\in(0,1)$. Then, $\varphi(x)$ decreases on $(x_q,+\infty)$ for all $q\in(0,1)$, and the conclusion follows.

3) For $q=1$ the result remains true by using Lemma \ref{der}.


 4) For the monotony of the function $\varphi_1(x)$, we use the same method as above and Lemma 1 in \cite{alz5}.

  The result of Proposition \ref{prou} is not true for all $q>1$ and all $x>1$. But it remains true with some extra conditions on $q$. For sake of simplicity we prove it on separate proposition. To do this we need some auxiliary results.
  \begin{lemma}\label{lel} For all $q>0$, and all $k\geq 1$, the function $\frac{\psi^{(k+1)}_q(x)}{\psi^{(k)}_q(x)}$ increases on $(0,+\infty)$.
  \end{lemma}
  {\bf Proof.} To prove the lemma, it sufficient to show that $S_{q,k}(x):=\psi^{(k+2)}_q(x)\psi^{(k)}_q(x)-(\psi^{(k+1)}_q(x))^2>0$ on $(0,+\infty)$ for all $q>0$. By the series expansion of the $q$-polygamma function, we have for $q\neq 1$, and $k\geq 2$
  $$S_{q,k}(x)=(\log q)^{2k+4}\Big(\sum_{n,m=1}^\infty n^{k+1}m^k(n-m)\frac{q^{(n+m)x}}{(1-q^n)(1-q^m)}\Big).$$
  Hence,
  $$S_{q,k}(x)=(\log q)^{2k+4}\Big(\sum_{n<m}^\infty n^km^k(n-m)^2\frac{q^{(n+m)x}}{(1-q^n)(1-q^m)}\Big)>0.$$
  For $k=1$, and $q\in(0,1)$, the proof reminds the same as the above case.

  For $k=1$ and $q>1$, we have $$S_{q,1}(x)=\psi'''_q(x)\psi'_q(x)-(\psi''_q(x))^2,$$
  hence, $$S_{q,1}(x)=\psi'''_{\frac1q}(x)(\psi'_{\frac1q}(x)+\log q)-(\psi''_{\frac1q}(x))^2=S_{\frac1q,1}(x)+\psi'''_{\frac1q}(x)\log q\geq 0$$
  For $q=1$, one has by Lemma \ref{der}, $\lim_{q\to 1}\frac{\psi_q^{(k+1)}(x)}{\psi_q^{(k)}(x)}=\frac{\psi^{(k+1)}(x)}{\psi^{(k)}(x)}$, and the result follows from the first item.
 Which implies that the function $S_{q,k}(x)$ increases on $(0,+\infty)$ for all $q>0$ and all $k\geq 1$.

In the sequel we give some monotonic results involving the functions $q$-polygamma with respect to the variable $q$ and fixed $x$.
\begin{lemma}\label{loule} \ \begin{enumerate}
\item[$(1)$] For every $x>0$, the function $q\mapsto\psi'_q(x)$ increases on $(0,1)$.
\item[$(2)$] The function $q\mapsto\psi_q(x)-\psi_q(1)$ decreases on $(0,+\infty)$ for $x\in(0,1]$ and increases for $x\in[1,+\infty)$
\end{enumerate}
\end{lemma}
{\bf Proof.}
1) Recall that for $q\in(0,1)$ $$\psi_q(x)=-\log(1-q)+\sum_{n=0}^\infty\frac{(\log q)q^{n+x}}{1-q^{n+x}}.$$
Then, $$\psi'_q(x)=\sum_{n=0}^\infty\frac{(\log q)^2q^{n+x}}{(1-q^{n+x})^2}.$$
For $y>0$ and $q\in(0,1)$, let $\tau(q)=\frac{(\log q)^2q^{y}}{(1-q^{y})^2}$. Differentiation of $\tau(q)$ yields,
$$\tau'(q)= \frac{yq^{y-1}(1 + q^y) \log q}{(1 - q^y)^3} (\frac{2 - 2 q^y}{(1 + q^y) y } + \log q).$$
To show that $\tau(q)$ increases it is sufficient to show that $\beta(q)=\frac{2 - 2 q^y}{(1 + q^y) y } + \log q\leq 0$. Another calculation shows that $$\beta'(q)=\frac{( q^y-1)^2}{(q (1 + q^y)^2)}\geq0\quad q>0,$$
and $\beta(1)=0$, so $\beta(q)\leq 0$ for $0<q\leq 1$.  Thus $\tau(q)$ is an increasing function of $q$, and so $\psi'_q(x)$ is an increasing function of $q$, for $q\in(0,1)$ and all $x>0$.

2) Firstly, remark that the function $x\mapsto\psi_q(x)-\psi_q(1)$ is non-positive on $(0,1)$ and non-negative on $[1,+\infty)$.
Moreover,
$$\psi_q(x)-\psi_q(1)=\log q\big(x-1-\sum_{n=1}^\infty\frac{1-q^{n(1-x)}}{1-q^n}\big):=E(x,q)\log q.$$
Let $R(x,q)=\frac{1-q^{n(1-x)}}{1-q^n}$. Differentiate with respect to $q$ yields $$\partial_qR(x,q)=\frac{n q^{-1 + n - n x} (1-q^{n x} -x( 1 - q^n ))}{( q^n-1)^2}.$$
Since, $q\mapsto\partial_q (1-q^{n x} -x( 1 - q^n ))=nxq^{n-1}(1-q^{n(x-1)})\geq 0$ for $x\in(0,1]$ and $\leq 0$ for $x\geq 1$. Then, $1-q^{n x} -x( 1 - q^n )\geq 0$ for $x\in(0,1]$ and $\leq 0$ for $x\geq 1$ and all $q\geq 1$. Thus, the function $q\mapsto R(x,q)$ increases on $[1,+\infty)$ for all $x\in(0,1]$ and decreases for $x\geq 1$. Then, the function $q\mapsto E(x,q)$ decreases on $[1,+\infty)$ for all $x\in(0,1]$, and is non-positive. Moreover, it increases on $[1,+\infty)$ for all $x\geq 1$, and non-negative. Which implies that $q\mapsto \psi_q(x)-\psi_q(1)$ decreases on $[1,+\infty)$ for $x\in(0,1]$ and increases for all $x\geq 1$.

Let $q\in(0,1]$ and $x>0$, then $\displaystyle\psi_q(x)-\psi_q(1)=\int_1^x\psi_q'(t)dt$. Using item 1), one deduces that $\psi_q(x)-\psi_q(1)$ is an increasing function of $q$ for $x\geq 1$ and is a decreasing function of $q$ for $x\in(0,1]$.

\begin{lemma}\label{loul} \ \begin{enumerate}

\item [$(1)$] For $x>0$, the function $q\mapsto\psi''_q(x)$ decreases on $(0,1)$ and increases on $[1,+\infty)$.
\item [$(2)$] For $x>0$, The function  $q\mapsto\psi'_q(x)$ increases on $(0,+\infty)$.
\item[$(3)$] The function $q\mapsto\psi_q(x)$ decreases on $(0,\infty)$ for all $x\in(0,1]$ and increases on $(0,+\infty)$ for all $x\geq 2$.
\end{enumerate}
\end{lemma}
{\bf Proof.}
1) For $q\in(0,1)$, and $x>0$, we have, $$\psi_q''(x)=\sum_{n=0}^\infty\frac{q^{n+x}(1+q^{x+n})(\log q)^3}{(1-q^{n+x})^3}.$$
For $q\in(0,1)$ and $a>0$, let $$h(a,q)=\frac{q^{a}(1+q^{a})(\log q)^3}{(1-q^a)^3},$$
A first differentiation with respect to $q$ gives $$\begin{aligned}\partial_qh(q,a)&=\frac{q^{-1 + a}
  (\log q)^2 a (1 + q^a (4 + q^a))}{(-1 + q^
  a)^4} (\frac{3 - 3 q^{2 a}}{ a (1 + q^a (4 + q^a))} + \log q))\\&=\frac{q^{-1 + a}
  (\log q)^2 a (1 + q^a (4 + q^a))}{(-1 + q^
  a)^4} m(q,a).\end{aligned}$$
  Furthermore, $$\partial_q m(q,a)=\frac{(-1 + q^a)^4}{q (1 + q^a (4 + q^a))^2},$$ and $m(1,a)=0$ for all $a>0$. Hence, $\partial_qh(q,a)\geq 0$ for $q\geq 1$  and $\partial_qh(q,a)\leq 0$ for $q\in(0,1]$ and for all $a>0$. Thus, $q\mapsto h(q,a)$ decreases on $(0,1]$ and increases on $[1,+\infty)$ for all $a>0$.
 Which gives the desired result.

2) Let $0<p<q$. By the mean value theorem we have for all $x>0$ there is $c$ such that
$$\psi'_p(x)-\psi'_q(x)= \psi'_p(1)-\psi'_q(1)+(x-1)(\psi''_p(c)-\psi''_q(c)).$$

a) If $x\in (0,1)$ then $c\in (x,1)$. Moreover, by the previous item $\psi''_p(c)-\psi''_q(c)\geq 0$ and $\psi'_p(1)-\psi'_q(1)\leq 0$. Then $\psi'_p(x)-\psi'_q(x)\leq 0$.

b) If $x\in (1,+\infty)$ then $c\in (1,x)$. Moreover, by the previous item $\psi''_p(c)-\psi''_q(c)\leq 0$ and $\psi'_p(1)-\psi'_q(1)\leq 0$. Again $\psi'_p(x)-\psi'_q(x)\leq 0$.

3) We saw that $q\mapsto\psi_q(1)$ decreases and by item 2) of Lemma \ref{loule}, the function $q\mapsto\psi_q(x)-\psi_q(1)$ decreases on $(0,1]$ for all $x\in(0,1]$. Then, $q\mapsto\psi_q(x)$ is a decreasing function of $q$, $q\in(0,1]$ and all $x\in(0,1]$.

Furthermore, for every $x\geq 2$, $\psi_q(x)=\psi_q(2)+\int_2^x\psi'_q(t)dt$. Since, the function $q\mapsto\psi_q(2)$ increases, using item 2), we deduce the desired result.
\begin{corollary}\label{loul} \
 For $x\geq 2$, the function  $q\mapsto x\psi'_q(x)+2\psi_q(x)$ increases on $(0,\infty)$.
In particular, for any $q \geq1$ and $x \geq 2$ we have $$x\psi'_q(x)+2\psi_q(x)\geq x\psi'(x)+2\psi(x).$$
\end{corollary}
{\bf Proof.}
1) Let $u(q,x)=x\psi_q'(x)+2\psi_q(x)$. Differentiate $u(q,x)$ with respect to $x$ yields
$\partial_x u(q,x)=3\psi'_q(x)+x\psi''_q(x).$ In Lemma \ref{loul} it is proved that $q\mapsto\psi''_q(x)$ increases on $[1,+\infty)$ and decreases on $(0,1)$, $q\mapsto\psi'_q(x)$ increases on $(0,+\infty)$ for $x>0$ and the function $x\mapsto\psi_q(x)$ increases on $(0,+\infty)$ for $x\geq 2$. Then, for $1\leq p<q$ we have
$\partial_x u(p,x)\leq \partial_x u(q,x)$. Integrate on $[2,x]$ gives
$$u(p,x)-u(q,x)\leq u(p,2)-u(q,2).$$
Which is non-positive. Then, for every $x\geq 2$, $u(p,x)\leq u(q,x)$.




 \begin{proposition}\label{proi}\ Let $q>0$. \begin{enumerate}
  \item The function $G_q(x)=\frac{\psi'_q(x)}{\psi_q(x)}$ decreases on $(1,+\infty)$ for all $q\geq1$.

 \item  The function $\varphi_q(x)=x\frac{\psi'_q(x)}{\psi_q(x)}$ decreases on $(1,x_q)$ for all $q>1$, and it decreases on $(1,+\infty)$ if and only if $q\in(0,1]\cup[q_0,+\infty)$, where \\$q_0=\frac1{3\sqrt[3]2}(2\sqrt[3]2 + \sqrt[3]{25-3\sqrt{69}}+\sqrt[3]{25+3\sqrt{69}})\simeq 1.75488$ is the unique positive solution of $(\sqrt q)^3-\sqrt q-1=0$.
  \end{enumerate}
  \end{proposition}
  {\bf Proof.} 1)
Differentiation of $G_q(x)$ gives
$$G'_q(x)=\frac{\psi''_q(x)\psi_q(x)-(\psi'_q(x))^2}{(\psi_q(x))^2},$$

Since, $\psi''_q(x)\leq 0$ and for $x>x_q$, $\psi_q(x)> 0$ for all $q>0$. Hence, $G'_q(x)\leq 0$ and $G_q(x)$ decreases on $(x_q,+\infty)$ for all $q>0$

Now assume that $q>1$ and $x\in[1,3/2]$, $x\neq x_q$. On the first hand we have
$$G'_q(x)=\frac{\psi''_{\frac1q}(x)\big((x-\frac32)\log q+\psi_{\frac1q}(x)\big)-(\psi'_{\frac1q}(x)+\log q)^2}{(\psi_q(x))^2}.$$
Easy computation yields,
$$G'_q(x)=\frac{\psi''_{\frac1q}(x)\psi_{\frac1q}(x)-(\psi'_{\frac1q}(x))^2}{(\psi_q(x))^2}
+H_q(x),$$
where $$H_q(x)=\frac{\psi''_{\frac1q}(x)(x-\frac32)\log q-2(\log q)\psi'_{\frac1q}(x)-(\log q)^2}{(\psi_q(x))^2}.$$
 Using Proposition \ref{prou}, we get $G'_q(x)\leq H_q(x)$.

 For $\alpha\in(0,1)$, and $x\in (1,\frac32)$, setting, $$K_\alpha(x)=\psi''_{\alpha}(x)(x-\frac32)-2\psi'_{\alpha}(x)+\log \alpha.$$
Easy computation gives,
$$K_\alpha(x)=(\log\alpha)^2\sum_{n=0}^\infty\big((x-\frac32)\frac{1+\alpha^{x+n}}{1-\alpha^{n+x}}\log\alpha-2\big)\frac{\alpha^{x+n}}{1-\alpha^{n+x}}
+\log\alpha.$$
For $x\in[1,3/2]$ and $\alpha\in(0,1)$, let $$\kappa_\alpha(x)=\big((x-\frac32)\frac{1+\alpha^{x+n}}{1-\alpha^{n+x}}\log\alpha-2\big)\frac{\alpha^{x+n}}{1-\alpha^{n+x}}.$$
We have, $$\kappa'_\alpha(x)=\frac{\alpha^{n + x}
  \log\alpha \Big(2 (-1 + \alpha^{2 (n + x)}) + (1 +
      \alpha^{n + x} (4 + \alpha^{n + x})) (-3 + 2 x) \log\alpha\Big)}{2 (-1 + \alpha^{
   n + x})^4}:=\frac{(\alpha^{n + x}
  \log\alpha)\mu_\alpha(x) }{2 (-1 + \alpha^{
   n + x})^4},$$
   and $$\mu'_\alpha(x)=\log\alpha \big(2 + 8 \alpha^{n + x} + 6 \alpha^{2 (n + x)} +
   2 \alpha^{n + x} (2 + \alpha^{n + x}) (-3 + 2 x) \log\alpha\big).$$
   It seems to be clear that $\mu'_\alpha(x)<0$ for $x\in[1,3/2]$ and $\alpha\in(0,1)$, thus $\mu_\alpha(x)$ decreases on $[1,3/2]$. Furhermore,
   $\mu_\alpha(3/2)<0$ and $\mu_\alpha(1)=2 (-1 + \alpha^{2 (n + 1)}) - (1 +
      \alpha^{n + 1} (4 + \alpha^{n + 1}))  \log\alpha.$

If $\mu_\alpha(1)\leq 0$, then $\kappa'_\alpha(x)\geq 0$ for all $x\in[1,3/2]$ and then, $\kappa_\alpha(x)\leq \kappa_\alpha(3/2)\leq 0$ and $K_\alpha(x)\leq 0$.

If $\mu_\alpha(1)>0$, then there is a unique $x_0\in(1,3/2)$ such that $\kappa_\alpha(x)$ decreases on $(1,x_0)$ and increases on $(x_0,3/2)$. That is to say, $\kappa_\alpha(x)\leq\max(\kappa_\alpha(1),\kappa_\alpha(3/2))$. One can see that, $\kappa_\alpha(3/2)\leq 0$, thus, $K_\alpha(x)\leq\max(0, K_\alpha(1))$. So, to prove that $K_\alpha(x)\leq0$. It sufficiencies to prove that  $K_\alpha(1)\leq 0$ for all $\alpha\in(0,1)$.

We have, $$K_\alpha(1)=-\frac12(\psi''_\alpha(1)+4\psi'_\alpha(1)-2\log\alpha).$$
Since, $\psi'_\alpha(1)\geq 0$ and $\psi'_\alpha(1)=\psi'_\frac1\alpha+\log\alpha$, then $K_\alpha(1)\leq -\frac12(\psi''_\alpha(1)+2\psi'_{\frac1\alpha}(1)).$ Using Lemma \ref{loul} and Lemma \ref{der}, then the function $\alpha\mapsto \psi''_\alpha(1)+2\psi'_{\frac1\alpha}(1)$ decreases on $(0,1)$, and
$$\psi''_\alpha(1)+2\psi'_{\frac1\alpha}(1)\geq \psi''(1)+2\psi'(1)>0.$$

2) a) Assume $q>1$ and $x\in(1,x_q)$. Since, $-\varphi_q(x)=x(-\frac{\psi'_q(x)}{\psi_q(x)})$, hence $-\varphi_q(x)$ is a product of two positive increasing functions. Then $\varphi_q(x)$ decreases on $(1,x_q)$.

b) For $x> x_q$, $$\varphi_q(x)=\frac{x\psi'_{\frac 1q}(x)}{\psi_q(x)}+\log q\frac{x}{\psi_q(x)}.$$
Since, by Lemma \ref{l3} the function $x\psi'_{\frac 1q}(x)$ decreases and is positive and $\frac1{\psi_q(x)}$ decreases and is positive, then $\frac{x\psi'_{\frac 1q}(x)}{\psi_q(x)}$ decreases on $(x_q,+\infty)$.

Let $v_q(x)=\frac{x}{\psi_q(x)}$, then $(\psi_q(x))^2v_q'(x)=\psi_q(x)-x\psi_q'(x)$. Moreover, the derivative of the right hand side is $-x\psi''_q(x)$ which is positive. Then the function $\psi_q(x)-x\psi'_q(x)$ increases on $(x_q,+\infty)$, furthermore, by easy computation we have $\displaystyle \lim_{x\to+\infty}\psi_q(x)-x\psi'_q(x)=-\log(\sqrt q(q-1))\leq0$ for $q\geq q_0$. The case $q\in(0,1)$ follows from Proposition \ref{prou}. One deduces that $\varphi_q(x)$ decreases.

{\it The converse.} Assume that $\varphi_q(x)$ decreases.
Then, $$(\psi_q(x))^2\varphi'_q(x)=\psi'_q(x)(\psi_q(x)-x\psi'_q(x))+x\psi_q(x)\psi''_q(x)\leq 0.$$
For $q>1$, and $x>1$, we have
$$\sum_{n=1}^\infty n^2\frac{q^{-nx}}{1-q^{-n}}\leq q^{-x}\frac q{q-1}(1+(q-1)\sum_{n=2}^\infty n^2\frac{1}{1-q^{-n}}).$$
Then, \begin{equation}\label{p}|\psi''_q(x)|\leq a_q q^{-x}.\end{equation}
Following the same method, there is $b_q,c_q\geq 0$, such that,
\begin{equation}\label{q}|\psi_q(x)|\leq b_q+(\log q)x+c_q q^{-x}.\end{equation}
From inequality \eqref{p} and \eqref{q}, we get $\displaystyle\lim_{x\to\infty}x\psi_q(x)\psi_q''(x)=0$, and $\lim_{x\to+\infty}\frac{x\psi_q''(x)}{\psi_q(x)}=0$. Also we saw that $\displaystyle\lim_{x\to+\infty}\psi_q'(x)=\log q$, then
$$-(\log q)\log(\sqrt q(q-1))\leq 0,$$
or equivalently
$$q^{\frac32}-q^{\frac12}-1\geq 0.$$
Which implies that $q\in[q_0,+\infty)$.
\begin{figure}[h]
\centering\scalebox{0.5}{\includegraphics[width=17cm, height=12cm]{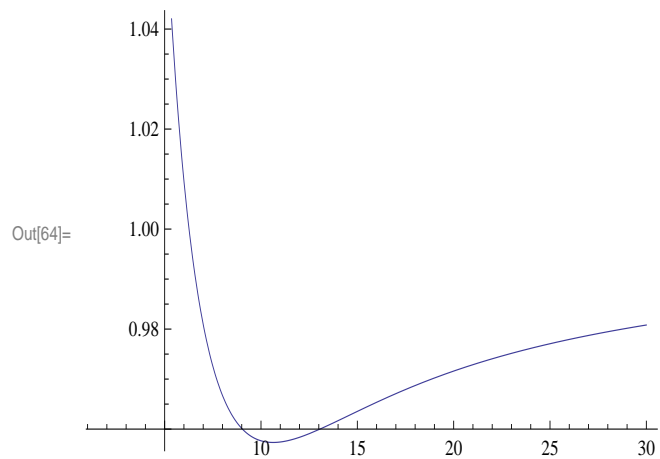}}
\caption{$\varphi_{1.6}(x)=\frac{x\psi'_{1.6}(x)}{\psi_{1.6}(x)}$}
\end{figure}
\newpage
  \begin{corollary} For every $q>0$, the function $|\psi_q(x)|$ is logarithmic concave on $(1,+\infty)$.
  \end{corollary}
  {\bf Proof.} For $x>0$  $x\neq x_q$, let $h_q(x)=\log(|\psi_q(x)|)$. Then, $h'_q(x)=G_q(x)$, and $h_q''(x)=G'_q(x)\leq 0$, and the result follows by Proposition \ref{prou} and  Proposition \ref{proi}.

In the proposition below we provided an extension of the result of Lemma \ref{i}.
\begin{proposition} For $x\geq 1$, $q>0$ and $a\in\Bbb R$, we set $h(x,a,q)=x\psi'_q(x)+a\psi_q(x)$. Then
For a given $q>0$, $h(x,a,q)\geq 0$ for all $x\geq 1$ if and only if $0\leq a\leq -\frac{\psi'_q(1)}{\psi_q(1)}$.

In particular $\psi'_q(x)+\psi_q(x)\geq 0$ and $\psi'_q(x)+2\psi_q(x)\geq 0$ for all $q>0$.
\end{proposition}
The proposition gives a refinement of the result of Alzer Lemma \ref{i}. Indeed, $-\frac{\psi_q'(1)}{\psi_q(1)}\geq 2$ for all $q>0$.

{\bf Proof.} By Proposition \ref{proi} and \ref{prou}, we have the function $\psi'_q(x)/\psi_q(x)$ decreases on $[1,+\infty)$ for all $q>0$. Then, for $x\in[1,x_q)$,
$$\frac{\psi'_q(x)}{\psi_q(x)}-\frac{\psi'_q(1)}{\psi_q(1)}\leq 0.$$
So, for $x\in[1,x_q)$, $$\psi'_q(x)+a\psi_q(x)=\psi_q(x)(\frac{\psi'_q(x)}{\psi_q(x)}+a)\geq \psi_q(x)(\frac{\psi'_q(1)}{\psi_q(1)}+a)\geq 0.$$
Thus, for $x\geq x_q$, $\psi'_q(x)\geq 0$ and $a\psi_q(x)\geq 0$ and the result follows.

For the converse. If for all $x\geq 1$, $\psi'_q(x)+a\psi_q(x)\geq 0,$ then for $x=1$, we get $a\leq -\frac{\psi'_q(1)}{\psi_q(1)}$. Moreover, as $x\to+\infty$, we get for $q\in(0,1)$, $-a\log(1-q)\geq 0$ and then $a\geq 0$. On the other hand, for $q>1$, we have, $\displaystyle\lim_{x\to\infty}\psi'_q(x)/\psi_q(x)=0$, then $a\geq 0$.

Remark that $\psi'_q(1)+2\psi_q(1)\geq 0$, then $-\psi'_q(1)/\psi_q(1)\geq 2$. On deduces the result for $a=1$ and $a=2$.


\section{Harmonic mean of the $q-$gamma function}
Our first main result is a generalization of Gautschi inequality in \cite{gaut}. We start with some useful lemmas.

Let's $$J=\Big\{q>0;\psi_q(1)-(\psi_q(1))^2+\psi'_q(1)\geq0\Big\}.$$

By Lemma \ref{l1}, we have $\psi'_q(1)\geq -2\psi_q(1)$, then $\psi_q(1)-(\psi_q(1))^2+\psi'_q(1)\geq -\psi_q(1)(1+\psi_q(1))$. Furthermore, one shows by induction that $4^n-1\geq (9/10) 4^n$ for $n\geq 2$. Then,
$$\psi_4(1)\geq \log \frac23-\frac23\log 2-\frac{20}9\log 2\sum_{n=2}^\infty\frac1{4^n}.$$
Thus, $$1+\psi_4(1)\geq 1+\log \frac23-\frac23\log 2-\frac{5}{27}\log 2\simeq0.00407>0.$$
Since, $\psi_q(1)<0$ and the function $q\mapsto1+\psi_q(1)$ decreases on $(0,+\infty)$, then $[0,4]\subset J$

Numerical computation show that $\psi_{10}(1)-(\psi_{10}(1))^2+\psi'_{10}(1)\simeq-0.072$, then $J \varsubsetneq [0,10)$.
\begin{lemma}\label{leu} For $q>0$, and $x\geq 1$, let $\theta_1(x)=\frac{x\psi_q(x)}{\Gamma_q(x)}$. Then
\begin{enumerate}
\item[$1)$] for $q\in J$, $\theta_1(x)$ increases on $[1,x_q]$.
\item[$2)$] for $q\notin J$ there is a unique $y_q\in(1,x_q)$ such that, $\theta_1(x)$ decreases on $(1,y_q)$ and increases on $[y_q,x_q]$.
\end{enumerate}
\end{lemma}
{\bf Proof.} Differentiation of $\theta_1(x)$ gives,
$$\theta'_1(x)=\big(x\psi_q(x)+x^2\psi_q'(x)-(x\psi_q(x))^2\big)\frac{1}{x\Gamma_q(x)}.$$
1) Let $q\geq 1$,  by Lemma \ref{l2} the function $x\psi_q(x)$ increases, moreover, it is non positive on $[1,x_q]$ then $x^2(\psi_q(x))^2$ decreases. By Lemma \ref{l3}, the function $x^2\psi'_q(x)$ increases. Then, for all $q>1$ and all $x\in [1,x_q]$
$$x\psi_q(x)+x^2\psi_q'(x)-(x\psi_q(x))^2\geq \psi_q(1)-(\psi_q(1))^2+\psi'_q(1).$$
The right hand side is positive for every $q\in J\cap[1,+\infty)$.
We conclude that that $\theta_1(x)$ increases on $(1,x_q)$ for all $q\in J\cap[1,+\infty)$.

If $q\in(0,1]$, we write $$\theta'_1(x)=\big(\psi_q(x)+x\psi_q'(x)-x(\psi_q(x))^2\big)\frac{1}{\Gamma_q(x)},$$
and we use the inequality $x\psi_q'(x)+2\psi_q(x)\geq 0$ to get
$$\theta'_1(x)=\big(1+x\psi_q(x)\big)\frac{-\psi_q(x)}{\Gamma_q(x)}.$$
Since, $x\psi_q(x)$ increases on $(1,+\infty)$, then $1+x\psi_q(x)\geq 1+\psi_q(1)\geq 1-\gamma>0$, and for $x\in(1,x_q]$ $\psi_q(x)\leq 0$. Which implies that $\theta_1(x)$ increases on $(1,x_q)$ for all $q\in(0,1]$.

2) If $q\notin J$, then $q>1$ and $\psi_q(1)-(\psi_q(1))^2+\psi'_q(1)<0$. Moreover, the function $x\psi_q(x)+x^2\psi_q'(x)-(x\psi_q(x))^2$ increases on $[1,x_q]$ and is positive at $x=x_q$, then there is a unique $y_q\in(1,x_q)$ such that $\theta_1(x)$ decreases on $(1,y_q)$ and increases on $(y_q,x_q)$.

 \begin{proposition} \ \label{p2}
 For $q>0$, $x>0$ and $\alpha>0$, let
 $$f_q(x)=\frac{\Gamma_q(x)\Gamma_q(1/x)}{\Gamma_q(x)+\Gamma_q(1/x)},$$
 and $$g_{q,x}(\alpha)=f_q(x^\alpha).$$
 \begin{enumerate} \item[$1)$]
 \begin{enumerate}
 \item[$(a)$] For $q\in J$,
 The function $f_q(x)$
decreases on $(0,1]$ and increases on $[1,+\infty)$.
 \item[$(b)$] For $q\notin J$,
 the function $f_q(x)$
decreases on $(0,1/{y_q}]\cup[1,y_q]$ and increases on $[1/{y_q}, 1]\cup[y_q,+\infty)$.

\end{enumerate}
\item[$2)$]
\item[$(a)$] For every $x> 0$ and $q\in J$, the function $g_{q,x}(\alpha)$ increases on $(0,+\infty)$.
\item[$(b)$] For every $x> 0$ and $q\notin J$ the function $g_{q,x}(\alpha)$ decreases on $(\frac1{y_q},y_q)$ and increases on $(0,\frac1{y_q}]\cup[y_q,+\infty)$.

 In particular,
 For every $q\in J$, and $x>0$, $$2\frac{\Gamma_q(x)\Gamma_q(1/x)}{\Gamma_q(x)+\Gamma_q(1/x)}> 1.$$
 The sign of equalities hold if and only if $x=1$.

 For $q\notin J$, and $x>0$,
 $$\frac{\Gamma_q(x)\Gamma_q(1/x)}{\Gamma_q(x)+\Gamma_q(1/x)}\geq f_q(y_q).$$
 \end{enumerate}
\end{proposition}
One shows that $y_q\to 1$ as $q\to 1$.

{\bf Proof.} 1) A direct calculation gives
 $$f_q'(x)=(\theta_1(x)-\theta_1(\frac1x))\frac{f_q(x)}x,$$
where $\theta_1(x)=x\frac{\psi_q(x)}{\Gamma_q(x)}.$

1) Let $q\in J$. Since, for $x\geq 1$, we have $x\geq \frac1x$ and by performing the relation between $f_q$ and $\theta_1$ and Lemma \ref{leu}, we get $f_q'(x)\geq 0$ for $x\in[1,x_q)$. Now, for $x>x_q$,  $\theta_1(x)\geq 0$, $\theta_1(1/x)\leq 0$ and  $f_q(x)\geq 0$. Then $f_q$ increases on $[1,+\infty)$. Moreover, $f_q(x)=f_q(1/x)$, hence $f_q(x)$ decreases on $(0,1]$.

For $q\notin J$. Applying Lemma \ref{leu}, we get for $x\in(1,y_q)$, $f'_q(x)\leq 0$ and $f'_q(x)\geq 0$ on $(y_q,x_q)$. It follows that $f_q(x)$ decreases on $(1,y_q)$ and increases on $[y_q,x_q]$. If $x\geq x_q$, as above we have $f'_q(x)\geq 0$ and $f_q(x)$ increases. By the relation $f_q(x)=f_q(1/x)$ we get the desired result.


2) Let $\varphi(\alpha)=\frac{\Gamma_q(x^\alpha)\Gamma_q(1/x^\alpha)}{\Gamma_q(x^\alpha)+\Gamma_q(1/x^\alpha)}$. On the first hand, we have
$$\varphi(\alpha)=f_q(x^\alpha),$$
and $\varphi'(\alpha)=x^\alpha\log(x)f_q'(x^\alpha).$
Applying item 1), we deduce that $\varphi'(\alpha)\geq 0$ for all $\alpha\geq 0$. This completes the proof.

As consequence we have the following two corollaries.
\begin{corollary}\label{p1} For every $q\in J$ and $x>0$, $$\Gamma_q(x)+\Gamma_q(\frac1x)\geq 2,$$
and $$\Gamma_q(x)\Gamma_q(\frac1x)\geq 1,$$
\end{corollary}

\begin{corollary} Letting $q\to 1$, we get that
$\displaystyle f(x)=\frac{\Gamma(x)\Gamma(1/x)}{\Gamma(x)+\Gamma(1/x)},$
decreases on $(0,1]$ and increases on $[1,+\infty)$.
\end{corollary}
\begin{figure}[h]
\centering\scalebox{0.5}{\includegraphics[width=17cm, height=12cm]{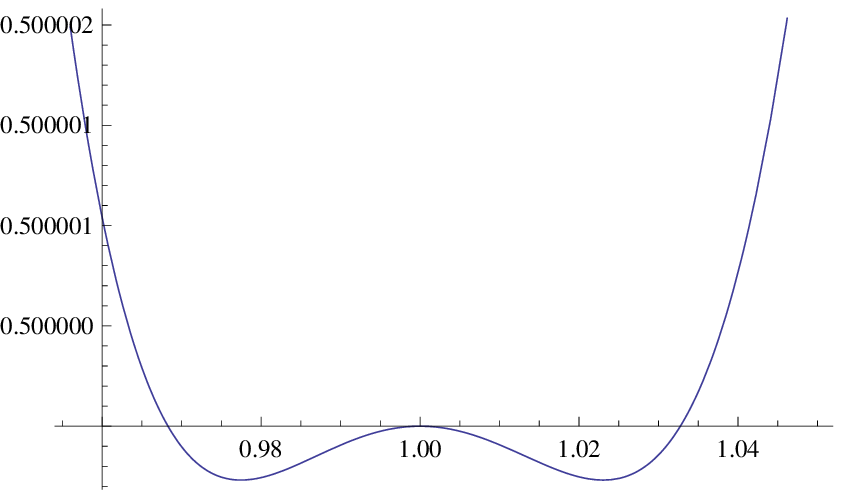}}
\caption{$f_{9.1}(x)$}
\end{figure}
\newpage
\begin{figure}[h]
\centering\scalebox{0.5}{\includegraphics[width=17cm, height=12cm]{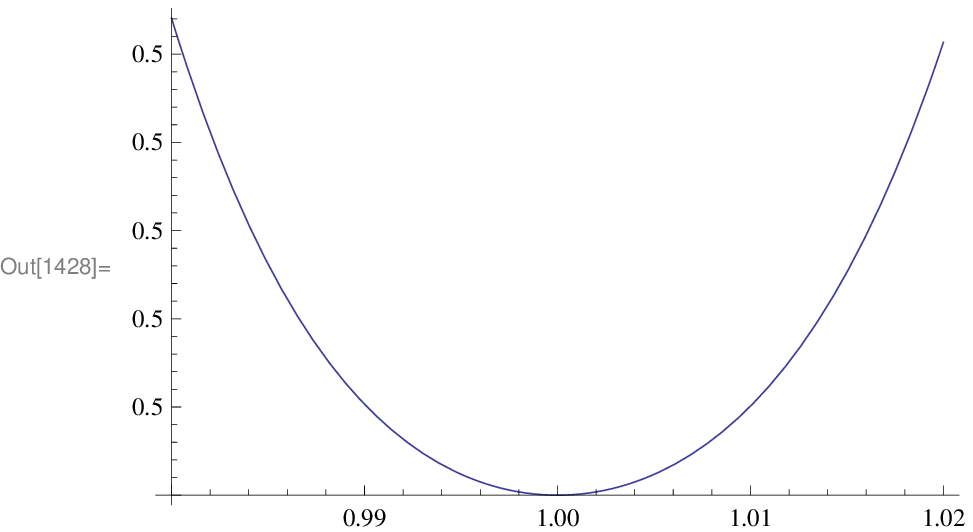}}
\caption{$f_{9}(x)$}
\end{figure}

Now we provide another generalization of the result of Proposition \ref{p2} when $q\in(0,1)$.

For $m\in\Bbb R$ and $a,b>0$, we set $H_m(a,b)=(\frac{a^m+b^m}{2})^{\frac1m}$.
\begin{proposition} \ Let $G_{m,q}(x)=H_m(\Gamma_q(x),\Gamma_q(\frac1x))$
\begin{enumerate}\item[$(1)$] For $q\geq1$,
the function $G_{m,q}(x)$ decreases on $(0,1)$ and increases on $(1,+\infty)$ if and only if $m\geq\frac{-\psi_q(1)-\psi'_q(1)}{(\psi_q(1))^2}$.
\item[$(2)$] For $q>0$ and $m\geq\frac1{\psi_q(1)}$ the function $G_{m,q}(x)$ decreases on $(0,1)$ and increases on $(1,+\infty)$.
\end{enumerate}
\end{proposition}
As a consequence, for $q=1$, we get $m\geq \frac1{\gamma}-\frac{\pi^2}{6\gamma}.$ This case is proved by Alzer \cite{alz4}. For $m=-1$, we recover the result of Proposition \ref{p2}.

Remark that $\psi_q'(1)+\psi_q(1)\geq \psi_q'(1)+2\psi_q(1)\geq 0$, then $-(\psi_q'(1)+\psi_q(1))/(\psi_q(1))^2\leq 0$ for all $q>0$.

{\bf Proof.} We use a method of Alzer in \cite{alz4}.
Let $ u_{m,q}(x)=\log G_{m,q}(x)$. If the inequality in the Proposition is true for all $x>0$, then, $u_{m,q}(x)\geq u_{m,q}(1)=0$, and $u'_{m,q}(1)=0$. Then,
$$u''_{m,q}(1)=(m-1)(\Gamma_q'(1))^2+\Gamma_q'(1)+\Gamma''_q(1)\geq 0.$$
From the relations, $\psi_q(1)=\Gamma'_q(1)$, and $\psi'_q(1)=\Gamma''_q(1)-(\Gamma'_q(1))^2$, we conclude that
$$m\geq -\frac{\psi_q(1)+\psi_q'(1)}{(\psi_q(1))^2}.$$

To prove the converse. First, remark that the function $m\mapsto H_m(a,b)$ increases on $\Bbb R$. So, it suffices to establish the inequality for $m=-\frac{\psi_q(1)+\psi_q'(1)}{(\psi_q(1))^2}.$ We shall prove that the function $G_{m,q}(x)$ satisfies, $G_{m,q}(x)\geq G_{m,q}(1)$. We have $g_{m,q}(x) = g_{m,q}(1/x)$, so that it is enough to show that $G_{m,q}(x)$ is strictly increasing
on $(1, +\infty)$.  Differentiation of $G_{m,q}(x)$ yields
$$G_{m,q}'(x)=\frac 12\Big((\Gamma_q(x))^{m-1}\Gamma'(x)-\frac1{x^2}(\Gamma_q(1/x))^{m-1}\Gamma'(1/x)\Big)(G_{m,q}(x))^{m-1}.$$
Which is equivalent to

$$G_{m,q}'(x)=\frac 1{2x}\Big(x(\Gamma_q(x))^{m}\psi_q(x)-\frac1{x}(\Gamma_q(1/x))^{m}\psi_q(1/x)\Big)(G_{m,q}(x))^{m-1}.$$
Let $v_{m,q}(x)=x(\Gamma_q(x))^{m}\psi_q(x)$. Then,

$$G_{m,q}'(x)=\frac 1{2x}\Big(v_{m,q}(x)-v_{m,q}(1/x)\Big)(G_{m,q}(x))^{m-1}.$$
If $x\geq x_q$, then $v_{m,q}(x)\geq 0$ and $v_{m,q}(1/x)<0$, hence $G'_{m,q}(x)>0$ and $G_{m,q}(x)$ increases on $[x_q,+\infty)$ for all $q>0$ and all $m\in\Bbb\R$.

Assume now $x\in(1,x_q)$. By Differentiation of $v_{m,q}(x)$ yields
$$v'_{m,q}(x)=(\Gamma_q(x))^{m}\big(\psi_q(x)+m x(\psi_q(x))^2+x\psi'_q(x)\big).$$

We use the inequality $x\psi'_q(x)+2\psi_q(x)\geq 0$ to get
$$v'_{m,q}(x)=-(\Gamma_q(x))^{m}\psi_q(x)\big(1-m x\psi_q(x)\big).$$
Since, for $m\geq 0$ and $x\in(1,x_q)$ the right hand side is positive for all $q>0$.

Let $m\in(-(\psi_q(1)+\psi'_q(1))/(\psi_q(1))^2, 0)$, By writing
$$v'_{m,q}(x)=\frac{(\Gamma_q(x))^{m}}x\big(x\psi_q(x)+m (x\psi_q(x))^2+x^2\psi'_q(x)\big),$$
and by Lemma \ref{l2} and Lemma \ref{l3}, we get for all $q>1$ and all $x\in(1,x_q)$
$$x\psi_q(x)+m (x\psi_q(x))^2+x^2\psi'_q(x)\geq \psi_q(1)+m (\psi_q(1))^2+\psi'_q(1)\geq 0.$$
 Which implies that $G_{m,q}(x)$ increases on $(1,+\infty)$.

2) Suppose now that $q>1$ and $m\geq1/\psi_q(1)$. Remark that
$$\frac1{\psi_q(1)}+\frac{\psi_q(1)+\psi'_q(1)}{(\psi_q(1))^2}=\frac{2\psi_q(1)+\psi'_q(1)}{(\psi_q(1))^2}\geq 0.$$
Then, $m\geq -(\psi_q(1)+\psi_q'(1))/(\psi_q(1))^2$ and the result is deduced by the first item.

Let $q\in(0,1)$, $m\geq1/\psi_q(1)$ and $x\in(1,x_q)$. Since,
$$v'_{m,q}(x)=(\Gamma_q(x))^{m}\big(\psi_q(x)+m x(\psi_q(x))^2+x\psi'_q(x)\big).$$

By the inequality $x\psi'_q(x)+2\psi_q(x)\geq 0$, we get

$$v'_{m,q}(x)\geq -\psi_q(x)(\Gamma_q(x))^{m}\big(1-m x\psi_q(x)\big).$$
If $m\geq 0$, then $1-m x\psi_q(x)\geq 0$ and $v'_{m,q}(x)\geq 0$.

If $m\in[1/\psi_q(1),0)$, by applying Lemma \ref{l2}, we get $1-mx\psi_q(x)\geq 1-m\psi_q(1)\geq 0$. Then
$$v_{m,q}'(x)\geq 0.$$

Then $G_{m,q}'(x)\geq 0$ and $G_{m,q}(x)$ increases on $(1,+\infty)$.

\section{Harmonic mean of the $q-$digamma function}
In this section we give some generalization of Alzer's and Jameson's inequalities proved in \cite{alz6}.
\begin{proposition}\label{p3} For all $x>0$ and $x\neq 1$, and $q\in(0,1)$ then
$$\psi_q(x)+\psi_q(\frac1x)<2\psi_q(1).$$
\end{proposition}
{\bf Proof.} Recall that, for $x>0$, $$\psi_q(x)=-\log(1-q)-\int_0^\infty\frac{e^{-xt}}{1-e^{-t}}d\gamma_q(t),$$
where $\gamma_q(t)=\left\{\begin{aligned}&-\log q\sum_{k=1}^\infty\delta(t+k\log q),\;0<q<1\\&t,\qquad q=1\end{aligned}\right.$

It is enough to prove the proposition for $x\in(0,1)$. Let
$$f(x)=\psi_q(x)+\psi_q(\frac1x).$$
Then, $$f'(x)=\psi'_q(x)-\frac1{x^2}\psi'_q(\frac1x)=\frac1x(u(x)-u(\frac1x)),$$
where $u(x)=x\psi'_q(x)$. Let $x\in(0,1)$, then $x<\frac1x$, hence by Lemma \ref{l3}
$$u(x)>u(\frac1x).$$
One deduces that, $f$ is strictly increasing on $(0,1)$. Then $f(x)<f(1)$ for all $x\in(0,1)$. By the symmetry $f(1/x)=f(x)$, one deduces the result on $(0,+\infty)$.
\begin{corollary} For all $x>0$, $x\neq 1$ and $q>1$,
$$\psi_q(x)+\psi_q(\frac1x)<\frac{(x-1)^2}x\log q+2\psi_q(1).$$
\end{corollary}
This corollary follows from Proposition \ref{p3} and the relation $\psi_q(x)=(x-3/2)\log q+\psi_{\frac1q}(x)$.

\begin{lemma}\label{pq} For every $x>0$ and $q>0$,
$$\psi_q'''(1+x)< -\frac{q^x(1+q^x)}{(1-q^x)^3}(\log q)^3<\psi_q'''(x),$$
$$\psi_q''(x)< -\frac{q^x}{(1-q^x)^2}(\log q)^2<\psi_q''(1+x).$$
\end{lemma}
{\bf Proof.} By the Lagrange mean value theorem there is $z\in(x,x+1)$ such that
$\psi''_q(1+x)-\psi''_q(x)=\psi'''_q(z).$ As the function $\psi'''_q(x)$ is strictly decreasing and by the relation
$\psi_q''(1+x)-\psi''_q(x)=-\frac{q^x(1+q^x)}{(1-q^x)^3}(\log q)^3$, one deduces the desired result.
The second inequality follows by a same method.

Let $$I=\{q>0,\psi'_q(1)+\psi''_q(1)\geq0 \}.$$
 \begin{lemma}\label{que} 
 There is a unique $p_0\in(1,\frac92)$,  such that $\psi'_q(1)+\psi''_q(1)\geq0$ if and only if $q\geq p_0$ and $I=[p_0,+\infty)$.

 Numerical computation shows that $p_0\simeq3.239945$.
\end{lemma}

{\bf Proof.} Let $$u(q):=\psi'_q(1)+\psi''_q(1),$$ and $I=\{q>0, u(q)\geq 0\}$. It was proved in Lemma \ref{l3} that $x\psi_q'(x)$ is strictly decreasing on $(0,+\infty)$ for all $q\in(0,1)$ then $\psi'_q(1)+\psi''_q(1)<0$. Furthermore, $\displaystyle\lim_{q\to 1}\psi'_q(1)+\psi''_q(1)=\psi'(1)+\psi''(1)=\zeta(2)-2\zeta(3)\simeq -0.759$. Then $I\subset (1,+\infty)$.

By Lemma \ref{loul} We saw that the function $q\mapsto\psi'_q(1)+\psi''_q(1)$ increases on $(1,+\infty)$.
Moreover, for $q\geq 2$, $$|\psi''_q(1)|=(\log q)^3\sum_{n=1}^\infty\frac{n^2}{q^n-1}\leq\frac{(\log q)^3}{q-1}\sum_{n=1}^\infty\frac{n^2}{2^n}.$$
Then, $\lim_{q\to+\infty}\psi''_q(1)=0$. Also,
$\psi'_q(1)\geq\log q,$ and $\lim_{q\to+\infty}\psi'_q(1)=+\infty$. Then, there is a unique $p_0>1$ such that $\psi'_q(1)+\psi''_q(1)< 0$ for $q\in(0,p_0)$ and $\psi'_q(1)+\psi''_q(1)\geq0$ for $p\geq p_0$.

By equation \eqref{i}, we get
$$\psi'_q(1)+\psi''_q(1)\geq \psi'_q(1)(1-\psi'_{\frac1q}(1)).$$
Since, $\psi'_q(1)\geq 0$ and the function $z(q)=1-\psi'_{\frac1q}(1)$ increases on $(0,+\infty)$. Furthermore,

$z(9/2)=1-(\log 9/2)^2\sum_{n=1}^\infty\frac n{(9/2)^n-1}$. Since, for $n\geq 2$, $(9/2)^n-1\geq  (6/5)4^{n}$. Then
$$z(9/2)\geq 1-\frac{(\log 9/2)^2}{7/2}-\frac56\frac{(\log 9/2)^2}{4}\sum_{n=2}^\infty\frac n{4^{n-1}},$$
Easy computation gives
$$z(9/2)\geq 1-\frac{(\log 9/2)^2}{7/2}-\frac{35(\log 9/2)^2}{216}\simeq0.067$$ Then,
$p_0<9/2$. This completes the proof.

Numerically
 $$u(3)\leq \log 3 + \frac12(\log 3)^2 (1 - \log 3) + \frac14(\log 3)^2 (1 - 2 \log 3) +
 \frac{27}{78} (\log 3)^2 (1 - 3 \log 3)\simeq-0.28132.$$
 Hence, $I\subset(3,9/2)$.
\begin{figure}[h]
\centering\scalebox{0.5}{\includegraphics[width=17cm, height=12cm]{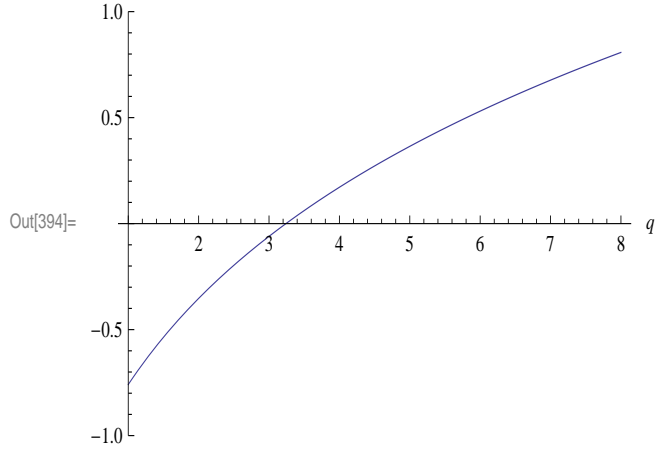}}
\caption{$\psi'_q(1)+\psi''_q(1)$}
\end{figure}

\newpage
\begin{lemma} \ \label{lo} \begin{enumerate}\item[$(1)$] For all $q>0$ and $x>0$, $2\psi''_q(x)+x\psi_q'''(x)\geq 0$.
\item[$(2)$] For $q\geq 1$, and $x>0$, $\psi'(x)+x\psi''(x)\leq\psi'_q(x)+x\psi_q''(x)\leq \log q$ and for $q\in(0,1)$ $\psi'_q(x)+x\psi_q''(x)\leq 0$
\item[$(3)$] The function $x\psi_q'(x)$ increases on $[1,+\infty)$ for every $q\in[p_0,+\infty)$ and decreases on $(0,1)$ if $q\in(0,p_0)$.
\end{enumerate}
\end{lemma}
{\bf Proof.}
1) Let $q\in(0,1)$ and $\varphi(x)=2\psi''_q(x)+x\psi_q'''(x)$, then
$$\varphi(1+x)-\varphi(x)=2\psi''_q(1+x)+(1+x)\psi_q'''(1+x)-2\psi''_q(x)-x\psi_q'''(x)$$
Since,
 $\psi_q'''(1+x)-\psi'''_q(x)=-\frac{q^x(1+q^x(4+q^x))}{(1-q^x)^4}(\log q)^4$, By Lemma \ref{pq} we have
$$\varphi(1+x)-\varphi(x)=\psi_q'''(1+x)-2\frac{q^x(1+q^x)}{(1-q^x)^3}(\log q)^3-x\frac{q^x(1+q^x(4+q^x))}{(1-q^x)^4}(\log q)^4.$$
By Lemma \ref{pq} we get
$$\varphi(1+x)-\varphi(x)\leq -\frac{q^x(\log q)^3}{(1-q^x)^4}(3(1+q^x)(1-q^x)+x\log q(1+q^x(4+q^x)).$$
For $u\in(0,1)$, let $j(u)=3(1+u)(1-u)+\log u(1+u(4+u))$. By successive differentiation we get
$j'(u)=4 + 1/u - 5 u + 2 (2 + u) \log u$, \\$j''(u)=-3 + (-1 + 4 u)/u^2 + 2 \log u$ and $j'''(u)=2 (-1 + u)^2/u^3>0$ on $(0,1)$. Then, $j''(u)\leq j''(1)=0$ and $j'(u)\geq j'(1)=0$. Thus, $j(u)$ increases on $(0,1)$ and $j(u)\leq j(1)=0$. Hence, for all $x>0$ and all $q\in(0,1)$
$$\varphi(1+x)-\varphi(x)\leq -\frac{q^x(\log q)^3}{(1-q^x)^4}j(q^x)\leq 0.$$
Or $\varphi(1+x)\leq\varphi(x)$ for all $x>0$ and then, for every $n\in\Bbb N$ and $x>0$
\begin{equation}\label{iqu}\varphi(x+n)\leq\varphi(x).\end{equation}
One shows by a similar method as in the Lemma \ref{pq} that $$\psi_q''(x)< -\frac{q^x}{(1-q^x)^2}(\log q)^2<\psi_q''(x+1).$$
From this relation together with Lemma \ref{pq}, one deduces that for $q\in(0,1)$ $\displaystyle\lim_{x\to+\infty}\psi''_q(x)=\displaystyle\lim_{x\to+\infty}x\psi'''_q(x)=0$.
By equation \ref{iqu} and as $n\to+\infty$ we get $\varphi(x)\geq 0$ for all $x>0$. Which gives the desired result.

For $q>1$, we saw that $2\psi''_q(x)+x\psi_q'''(x)=2\psi''_{\frac1q}(x)+x\psi_{\frac1q}'''(x)$ and the result follows.

2) Let $s(x)=x\psi_q'(x)$, then, $s'(x)=\psi'_q(x)+x\psi''_q(x)$ and $s''(x)=2\psi''_q(x)+x\psi'''_q(x).$
By the previous item we deduce that $s'(x)$ increases on $(0,+\infty)$ for all $q>0$. Since, $\displaystyle \lim_{x\to\infty}s'(x)=0$ if $q\in(0,1)$ and $=\log q$ if $q>1$. Which gives the desired result

3) We saw by item 2 that $s'(x)$ increases, then for every $x\geq 1$, $s'(x)\geq s'(1)=\psi'_q(1)+\psi''_q(1)\geq 0$ for all $q\in I$. Hence, $s(x)$ increases on $(1,+\infty)$ for $q\in I$.

\begin{proposition}\label{p05} For all $q\in [p_0,+\infty)$ and all $x>0$
$$\psi_q(x)+\psi_q(1/x)\geq 2\psi_q(1).$$
For $q\in(0,p_0)$,
$$\psi_q(x)+\psi_q(1/x)\leq 2\psi_q(1).$$
\end{proposition}
{\bf Proof.} Let $U(x)=\psi_q(x)+\psi_q(\frac1x)$, then $U'(x)=1/x(x\psi'_q(x)-1/x\psi'_q(1/x))$. If $x\geq 1$, then by Lemma \ref{lo}, and the fact that $x\geq 1/x$ we get $U'(x)\geq 0$ for all $q\geq p_0$. Hence, $U(x)$ increases on $(1,+\infty)$ and by the symmetry $U(x)=U(1/x)$, it decreases on $(0,1)$. Then $U(x)\geq U(1)$.

If $q\in(0,p_0)$, then $x\psi_q'(x)$ decreases on $(0,1)$, since $1/x\geq x$ then, $U'(x)\geq 0$ and $U(x)$ increases on $(0,1)$. By the symmetry $U(x)=U(1/x)$, it decreases on $(1,+\infty)$. Then $U(x)\geq U(1)$.Thus, $U(x)\leq U(1)$. Which completes the proof.
\begin{proposition}\label{p4} For all $x>0$ and $q>0$,
$$\psi_q(x)\psi_q(\frac1x)\leq(\psi_q(1))^2.$$
\end{proposition}

{\bf Proof.} Firstly, remark that the function $v(x)=\psi_q(x)\psi_{q}(\frac1x)$ is invariant by the symmetry $v(1/x)=v(x)$. So, it is enough to prove the result on $(1,+\infty)$ for all $q>0$.

By differentiation, we have
$$v'(x)=\frac1x(x\frac{\psi'_q(x)}{\psi_q(x)}-\frac1x\frac{\psi'_q(\frac1x)}{\psi_q(\frac1x)})v(x)=\frac1x(w(x)-w(\frac1x))v(x),$$
where $w(x)=x\frac{\psi'_q(x)}{\psi_q(x)}$.

  For $x\in(1,x_q)$, then $x>1/x$, and By Proposition \ref{prou} and Proposition \ref{proi} we have, the function $w(x)$ decreases. Then $w(x)\leq w(1/x)$. Moreover, $v(x)>0$, hence, $v'(x)<0$ and $v(x)<v(1)$.

For $x\geq x_q$, $\psi_q(x)\geq 0$ and $\psi_q(1/x)\leq 0$ hence, $v(x)\leq (\psi_q(1))^2$.

\begin{proposition}\label{p5} 1) For all $q\in (0,p_0)$  and all $x>0$
$$\frac{2\psi_q(x)\psi_q(\frac1x)}{\psi_q(x)+\psi_q(\frac1x)}>\psi_q(1).$$
2) For $q\in[p_0,+\infty)$, there is a unique $z_q>x_q$ such that for all $x\in[1/z_q,z_q]$
$$\frac{2\psi_q(x)\psi_q(\frac1x)}{\psi_q(x)+\psi_q(\frac1x)}>\psi_q(1).$$
If $x\in(0,1/z_q)\cup(z_q,+\infty)$, then the reversed inequality holds.

 The sign of equalities hold if and only if $x=1$.
\end{proposition}
{\bf Proof.} 1) Let $U(x)=\psi_q(x)+\psi_q(1/x)$. From Proposition \ref{p05}, we conclude that for $q\in(0,p_0)$ the expression $\frac1{U(x)}$ is
defined for all positive $x>0$. Applying Propositions \ref{p05} and \ref{p4}, we get for $q\in(0,p_0)$ and all $x>0$
$$\frac{2\psi_q(x)\psi_q(\frac1x)}{\psi_q(x)+\psi_q(\frac1x)}\geq \frac{2(\psi_q(1))^2}{\psi_q(x)+\psi_q(\frac1x)}>\psi_q(1).$$

2) From Proposition \ref{p05}, and the fact that for $q\geq p_0$, $U(x)$ increases on $(1,+\infty)$ and $U(1)=2\psi_q(1)<0$, $\lim_{x\to+\infty}U(x)=+\infty$, we deduce that there is a unique $z_q\in (1,+\infty)$ such that $U(z_q)=0$ and $U(x)$ is defined and negative for all $x\in(1/z_q,z_q)$. The fact that $z_q>x_q$ follows from the relation $\psi_q(z_q)=-1/\psi_q(1/z_q)>0$.

 Let $H(x)=\frac{\psi_q(x)\psi_q(\frac1x)}{\psi_q(x)+\psi_q(\frac1x)}$. Then,
$$H'(x)=x(\frac{x\psi'_q(x)}{(\psi_q(x))^2}-\frac{\psi'_q(1/x)}{x(\psi_q(1/x))^2})(H(x))^{-2}.$$
Since, $\psi_q(x)$ increases and is negative and by Proposition \ref{prou} and Proposition \ref{proi} $x\psi'_q(x)/\psi_q(x)$ decreases on $(1,x_q)$ for $q>1$ and is negative. Then, we get for $x>1$ $H'(x)>0$ and $H(x)$ increases on $(1,x_q)$. Thus, $H(x)\geq H(1)=1/2\psi_q(1)$.

If $x\in(x_q,z_q)$, then $\psi_q(x)\psi_q(1/x)<0$ and $\psi_q(x)+\psi_q(1/x)<0$, then $H(x)>0$. Hence, $H(x)>1/2\psi_q(1)$. One deduces the result on $(1/z_q,z_q)$ by the symmetry $H(x)=H(1/x)$.


If $x\in(0,1/z_q)\cup(z_q,+\infty)$, then $U(x)>0$. Moreover, $\psi_q(x)\psi_q(1/x)<0$. Then,
$\psi_q(x)\psi_q(1/x)\leq(\psi_q(1/x))^2$. Or equivalently
$$\frac{2\psi_q(x)\psi_q(\frac1x)}{\psi_q(x)+\psi_q(\frac1x)}\leq \psi_q(\frac1x).$$
Since, the function $x\mapsto\psi_q(1/x)$ decreases on $(0,+\infty)$ and $z_q>1$. One deduces the desired result.

As a consequence and since, $p_0>1$, by letting $q\to 1^-$, we get the following corollary
\begin{corollary} For all $x>0$,
$$\frac{2\psi(x)\psi(\frac1x)}{\psi(x)+\psi(\frac1x)}>-\gamma.$$
\end{corollary}

\end{document}